\documentclass{birkjour}
\usepackage{amsmath, amsfonts, amssymb, graphicx, color}
\usepackage{mathrsfs}
\usepackage{latexsym}

 \newtheorem{thm}{Theorem}[section]
 \newtheorem{cor}[thm]{Corollary}

 \theoremstyle{definition}
 \newtheorem{defn}[thm]{Definition}
 \theoremstyle{remark}

 \numberwithin{equation}{section}
\def\({\left ( }
\def\){\right )}
\def\<{\left < }
\def\>{\right >}

\begin{document}

\title[$\ast-\boldsymbol{\kappa}$-Ricci-Bourguignon.....]{Kenmotsu Contact Geometry Through the Lens of $\ast-\boldsymbol{\kappa}$-Ricci-Bourguignon Almost Solitons}

\author[K.Lavanya]{Lavanya Kumar}

\address{Department of Mathematics, School of Advanced Sciences, Vellore Institute of Technology, Chennai-600127, India}

\email{lavanya.k2023a@vitstudent.ac.in}

\author[S.Roy]{Soumendu Roy}

\address{Department of Mathematics, School of Advanced Sciences, Vellore Institute of Technology, Chennai-600127, India}

\email{soumendu.roy@vit.ac.in}




\begin{abstract}
 This paper focuses on the study of the newly introduced $\ast-\boldsymbol{\kappa}$-Ricci-Bourguignon almost soliton pertaining to  Kenmotsu structure manifolds. Our analysis concerns the characteristics of this soliton and derive the scalar curvature for a Kenmotsu manifold admitting such a structure. Further, we formulate the corresponding vector fields under the assumption that the manifold supports a $\ast-\boldsymbol{\kappa}-$Ricci-Bourguignon soliton. Additionally, we explore applications involving torse-forming vector fields within the framework of the $\ast-\boldsymbol{\kappa}-$Ricci-Bourguignon almost soliton on Kenmotsu structure manifolds. To support the theoretical findings, we provide a concrete illustration belonging to a $\ast-\boldsymbol{\kappa}-$Ricci-Bourguignon almost soliton in a 5D Kenmotsu structure manifold.
\end{abstract}

\maketitle

\section{Introduction}\label{1}
\hspace{.5cm}
 Contact geometry techniques play a major significance in the field of differential geometry and mathematical sciences. The fundamental of Classical mechanics present the development for the growth of Contact geometry. Ricci flow characterize the advancement of metrics described on connected almost contact metric manifolds that obtain automorphism groups of maximal dimension. Kenmotsu \cite{KK} introduced tensor equations characterizing a special class of manifolds in 1972. The manifolds were later known as Kenmotsu manifolds.\\
 
R. S. Hamilton \cite{HRS1} constructed the Ricci flow in 1982 as a differential equation determine the evolution of metrics on Riemannian manifolds. Suppose a compact Riemannian manifold $\mathbb{M}$ with metric $h$, in this case the ricci flow develops,
\begin{equation}\label{1.1}
 \frac{\partial h}{\partial \tau}= -2\mathbb{T}.
\end{equation}
A function to the Ricci flow (\cite{HRS1}, \cite{TP}) that develops only through a set of scaling and a one-parameter group of diffeomorphisms is known as a Ricci soliton \cite{HRS}. Its specifying equation is \cite{BRS}:\\
\begin{equation}\label{1.2}
 \mathscr{L}_\mathcal{V} h +2\mathbb{T}+2\mu h =0.
\end{equation}
In this condition, $\mathscr{L}_\mathcal{V} h$ represent to the Lie derivative of the metric $h$ respecting  the flow induced by the vector field $\mathcal{V}$, while $\mathbb{T}$ implies the tensor of Ricci curvature, $h$ acts as the  Riemannian metric, where $\mu$ stands for scalar constant.

The Ricci soliton is grouped into contracting, steady, or growing based on if the related scalar is negative, zero, or positive, correspondingly.\\

Wang, Gomes, and Xia \cite{XGW} have newly expanded the notion of a almost Ricci soliton by developing the general form of a almost Ricci soliton of type-$\boldsymbol{\kappa}$. A complete Riemannian geometric manifold equipped with $(\mathbb{M}^n, h)$ is known as almost Ricci soliton of type-$\boldsymbol{\kappa}$, written as $(\mathbb{M}^n, h, \mathcal{V}, \boldsymbol{\kappa}, \mu)$, admitting a smooth vector field $\mathcal{V}$ on $\mathbb{M}^n$, a function corresponding to the soliton $\mu \in C^\infty(\mathbb{M}^n)$, as well as a non-trivial real valued function $\boldsymbol{\kappa}$ on $\mathbb{M}^n$, where
    \begin{equation}\label{1.3}
        \boldsymbol{\kappa}\mathscr{L}_\mathcal{V} h +2\mathbb{T}+2\mu g =0.
    \end{equation}
    Similarly, to the notion of Ricci soliton \cite{HRS1}, and this signifies a self-repeating solution in the Ricci flow, there was a Ricci–Bourguignon soliton proposed through Shubham Dwivedi \cite{DS}:\\
         A Riemannian geometric manifold $(\mathbb{M}^n, h)$ can be classified as an RB soliton, namely the Ricci–Bourguignon soliton in the presence of an arbitrary vector field $\mathcal{V}$ associated with $\mathbb{M}$ that fulfills the condition,
      \begin{equation}\label{1.4}
      \mathscr{L}_\mathcal{V} h +2\mathbb{T}=2\Omega h+2\theta \mathcal{R} h,
      \end{equation}
      where $\mathscr{L}_\mathcal{V} h$ is the Lie derivative of the metric  $h$ concerning the vector field $\mathcal{V}$, along with $\Omega$ $\in$ $\mathcal{R}$ denotes a constant. 
       This is also stated as being contracting, steady, or growing with respect to $\Omega>0$, $\Omega=0$, together with $\Omega<0$, correspondingly.\\
       
Consider $(\mathbb{M},h)$ as a semi-Riemannian manifold \cite{CDM}.  A $\boldsymbol{\kappa}$-Ricci–Bourguignon almost soliton ($\boldsymbol{\kappa}$-RBA soliton) is characterized by the following equation:
\begin{equation}\label{1.5}
\boldsymbol{\kappa}\mathscr{L}_\mathcal{V} h+2\mathbb{T}+2\Omega h+2\theta \mathcal{R} h=0,
\end{equation}
  where $\mathcal{V}$ represents a smooth vector field, $\mathscr{L}_\mathcal{V} h$ symbolizes the Lie derivative of the metric $h$ relative to $\mathcal{V}$, $\mathbb{T}$ signifies the Ricci-type tensor of the manifold, $\mathcal{R}$ indicates the scalar curvature, $\boldsymbol{\kappa}$, $\lambda_1$, and $\beta_1$ are real constants.\\
 
The idea of the $\ast$-Ricci curvature tensor defined on almost Hermitian manifolds and the $\ast$-Ricci tensor corresponding to real hypersurfaces in non-flat complex space were presented by S.Tachibana \cite{TS} and T.Hamada \cite{TH}, correspondingly. The $\ast$-Ricci curvature tensor is specified as: 
\begin{equation}\label{1.6}
\mathbb{T}^\ast(\mathbf{P}, \mathbf{Q}) = \frac{1}{2}(Tr{\varphi \circ \mathcal{R}(\mathbf{P}, \varphi \mathbf{Q})})
\end{equation}
 for every vector field $\mathbf{P}$ and $\mathbf{Q}$ on $\mathbb{M}^n$, such that $Tr$ stands for trace and also $\varphi$ acts as a (1, 1)-tensor field. In the case that, $\mathbb{T}^\ast(\mathbf{P}, \mathbf{Q}) = \lambda h(\mathbf{P}, \mathbf{Q}) + \mathcal{V} \boldsymbol{\eta}(\mathbf{P}) \boldsymbol{\eta}(\mathbf{Q})$ for every pair of vector fields $\mathbf{P}$, $\mathbf{Q}$, with $\lambda$, $\mathcal{V}$ denoting smooth functions and $\boldsymbol{\eta}$ representing a linear-functional form(1-form), thereby, the manifold is called $\ast-\boldsymbol{\eta}$-Einstein manifold.\\
 Given that $\mathcal{V} = 0$, which means $\mathbb{T}^\ast(\mathbf{P}, \mathbf{Q}) = \lambda h(\mathbf{P}, \mathbf{Q})$ for arbitrary vector fields $\mathbf{P}$ as well as $\mathbf{Q}$, thus, the manifold is designated as $\ast$ - Einstein.\\
 In 2014, G. Kamis and K. Giotidou \cite{PKG} proposed the formulation of the $\ast$-Ricci soliton, described by the equation:
 \begin{equation}\label{1.7}
 \mathscr{L}_\mathcal{V} h + 2\mathbb{T}^\ast+ 2\mu h = 0 
 \end{equation}
 for any vector fields $\mathbf{P}$, $\mathbf{Q}$ on $\mathbb{M}^n$, with $\mu$ standing as a constant.\\
 Currently, Patra with co-authors. \cite{DPS} established the idea of almost $\ast$-Ricci–\\Bourguignon flow soliton($\ast$-RBS).\\
Every Riemannian geometric manifold $(\mathbb{M}^n, h)$ can be classified as $\ast$-Ricci–\\Bourguignon soliton or $\ast$-RB soliton, assuming the existence of a vector field $\mathcal{V}$ defined on $\mathbb{M}$ fulfilling the equation:
\begin{equation}\label{1.8}
     \mathscr{L}_\mathcal{V} h+ 2 \mathbb{T}^{\ast} = 2 [\Omega + \theta \mathcal{R} ^{\ast}] h,
 \end{equation}
 where $\mathscr{L}_\mathcal{V} h$ represents the Lie derivative of the metric $h$ concerning the vector function $\mathcal{V}$, $\Omega \in \mathcal{R} $ remains constant, $\mathbb{T}^\ast$ denotes the $\ast$-Ricci operator, and $\mathcal{R}^\ast=Tr(\mathbb{T}^\ast)$ signifies the $\ast$-scalar invariant of curvature.\\
 It is also described as contracting, steady or growing, corresponding to $\Omega \gtreqqless 0$, accordingly.\\
 If $\Omega$ serves as a continuous function, thus the $\ast$-RB soliton transforms into a almost $\ast$-RB soliton \cite{DPS}.\\
Using (1.3) together with (1.8), herein we define the $\ast$-$\boldsymbol{\kappa}$-Ricci-\\Bourguignon almost soliton  ($\ast$-$\boldsymbol{\kappa}RBS$) as follows:
\begin{defn}
 Any Riemannian geometric or, more generally, pseudo-\\Riemannian geometric manifold $(\mathbb{M}, h)$ of dimension $n$ is defined as $\ast$-$\boldsymbol{\kappa}$-Ricci-Bourguignon almost soliton ($\ast$-$\boldsymbol{\kappa}RBS$) provided that,
 \begin{equation}\label{1.9}
 \boldsymbol{\kappa}\mathscr{L}_\mathcal{V} h + 2\mathbb{T}^\ast-(2 \Omega + \theta \mathcal{R}^\ast)h = 0, 
 \end{equation}
 where $\mathscr{L}_\mathcal{V} h$ represents the Lie derivative of the metric $h$ following the vector field $\mathcal{V}$, $\mathbb{T}^\ast$ denotes $\ast$-Ricci tensor, $\mathcal{R}^\ast = Tr(\mathbb{T}^\ast)$ equals $\ast$-scalar curvature, together with $\mu, \alpha, \beta$ signify real scalars.\\
\end{defn}
 Conversely, a non-vanishing vector field $\rho$ on any Riemannian geometric or pseudo-Riemannian geometric manifold $(M,h)$ is termed Torse-forming tangent fields$(\rho)$ \cite{YK1} in case
 \begin{equation}\label{1.10}
 \nabla_\mathcal{V}\rho = \varphi \mathcal{V} + \theta(\mathcal{V})\rho, 
 \end{equation}
 where $\nabla$ denotes the Levi–Civita connection of $h$, $\varphi$ represents a smooth function, together with  $\theta$ represents a linear-functional form(1-form).  Furthermore, the field of vectors $\rho$ is termed concircular (\cite{BYC}, \cite{YK}) if the (1-form) $\theta$ is null everywhere within  (1.10). The field of vectors $\rho$ is termed concurrent (\cite{SJ}, \cite{CYK}) if, in (1.10), the (1-form) $\theta$ equals zero identically and where $\varphi=1$.  The vector quantity $\rho$ is known as recurrent the scalar function $\varphi = 0$ in (1.10).\\
 Whenever in (1.10) $\varphi= \theta$ equal zero, then the vector function $\rho$ is referred to as a uniform vector field .\\
 During 2017, B.Y.Chen \cite{BYC1} presented a novel vector quantity termed the torqued-generated vector field.  On condition that, the field of vectors $\rho$ satisfies (1.10) subject to $\theta(\rho)$vanishes, it follows that $\rho$ is designated torqued-generated vector field. For this case, $\rho$ regarded as torqued-related function, while the linear functional form(1-form) $\theta$ stands as torque-generated form of $\rho$.

\section{Preliminaries on Contact Metric Structures}\label{2}
\hspace{0.5cm}
Consider $\mathbb{M}$ as a connected manifold of dimension $(2n+1)$ endowed with an almost contact metric structure characterized by an almost contact Riemannian structure $(\varphi, \zeta, \boldsymbol{\eta}, h)$, $\varphi$ serves as a (1,1)-tensor field, $\zeta$ acts as a vector field, $\boldsymbol{\eta}$ stands for a linear-functional(1-form) form, and $h$ corresponds to the consistent Riemannian metric satisfying the following conditions: 
\begin{equation}
\varphi^2(\mathbf{P})=-\mathbf{P}+\boldsymbol{\eta}(\mathbf{P})\zeta, \boldsymbol{\eta}(\zeta)=1, \boldsymbol{\eta} \circ \varphi=0, \varphi \zeta=0,
\end{equation}
\begin{equation}
h(\varphi \mathbf{P}, \varphi \mathbf{Q})=h(\mathbf{P}, \mathbf{Q})-\boldsymbol{\eta}(\mathbf{P})\boldsymbol{\eta}(\mathbf{Q}),
\end{equation}
\begin{equation}
h(\mathbf{P}, \varphi \mathbf{Q})=-h(\varphi \mathbf{P}, \mathbf{Q}),
\end{equation}
\begin{equation}
h(\mathbf{P}, \zeta)=\boldsymbol{\eta}(\mathbf{P})
\end{equation}
for all vector fields $\mathbf{P}$, $\mathbf{Q}$ $\in \mathfrak{X}(\mathbb{M})$.\\
A almost contact-type Riemannian differentiable manifold can be described as a Kenmotsu structure on a manifold \cite{KK}. Supposing that the following conditions are satisfied:
\begin{equation}
(\nabla_\mathbf{P} \varphi)\mathbf{Q}=-h(\mathbf{P}, \varphi \mathbf{Q})\zeta-\boldsymbol{\eta}(\mathbf{Q})\varphi \mathbf{P} ,
\end{equation}
\begin{equation}
 \nabla_\mathbf{P} \zeta=\mathbf{P}-\boldsymbol{\eta}(\mathbf{P})\zeta ,
 \end{equation}
 Taking $\nabla$, the connection induced by the Riemannian metric $h$.\\
 Under the structure of a Kenmotsu-type manifold, the subsequent connections are valid (\cite{BPS}, \cite{BDR}):
 \begin{equation}
 \boldsymbol{\eta}(\mathcal{R}(\mathbf{P}, \mathbf{Q})\mathcal{R}) = h(\mathbf{P}, \mathcal{R})\boldsymbol{\eta}(\mathbf{Q}) - h(\mathbf{Q}, \mathcal{R})\boldsymbol{\eta}(\mathbf{P}),
 \end{equation}
 \begin{equation}
 \mathcal{R}(\mathbf{P}, \mathbf{Q})\zeta = \boldsymbol{\eta}(\mathbf{P})\mathbf{Q} - \boldsymbol{\eta}(\mathbf{Q})\mathbf{P},
 \end{equation}
 \begin{equation}
 \mathcal{R}(\mathbf{P}, \zeta)\mathbf{Q} = h(\mathbf{P}, \mathbf{Q})\zeta - \boldsymbol{\eta}(\mathbf{Q})\mathbf{P}, 
 \end{equation}
 where $\mathcal{R}$ denotes the Riemannian curvature tensor.
 \begin{equation}
 \mathbb{T}(\mathbf{P}, \zeta) = -2 n \boldsymbol{\eta}(\mathbf{P}),
 \end{equation}
 \begin{equation}
\mathbb{T}(\varphi \mathbf{P}, \varphi \mathbf{Q})
= \mathbb{T}(\mathbf{P}, \mathbf{Q})
+ 2 n \boldsymbol{\eta}(\mathbf{P})\, \boldsymbol{\eta}(\mathbf{Q}),
\end{equation}
 \begin{equation}
 (\nabla_\mathbf{P} \boldsymbol{\eta}) \mathbf{Q} = h(\mathbf{P}, \mathbf{Q}) - \boldsymbol{\eta}(\mathbf{P}) \boldsymbol{\eta}(\mathbf{Q})
 \end{equation}
 for any vector fields $\mathbf{P}, \mathbf{Q}, \mathbf{R} \in \mathfrak{X} (\mathbb{M})$.\\
 We have established that
 \begin{equation}
 (\mathscr{L}_\zeta h)(\mathbf{P}, \mathbf{Q}) = h(\nabla_\mathbf{P} \zeta, \mathbf{Q}) + h(\mathbf{P},\nabla_\mathbf{Q} \zeta)
 \end{equation}
 for any vector fields $\mathbf{P}, \mathbf{Q} \in \mathfrak{X} (\mathbb{M})$.\\
 Subsequently, we employ relations (2.6) and (2.13) to obtain
  \begin{equation}
 (\mathscr{L}_\zeta h)(\mathbf{P}, \mathbf{Q})=2[h(\mathbf{P}, \mathbf{Q})-\boldsymbol{\eta}(\mathbf{P})\boldsymbol{\eta}(\mathbf{Q})].
 \end{equation}
\textbf{Principle.2.1} (\cite{KNV}). In a specific (2n+1)-dim space Kenmotsu geometry, this particular $\ast$-Ricci operator represents expressed as
     \begin{equation}
 \mathbb{T}^\ast(\mathbf{P}, \mathbf{Q})=\mathbb{T}(\mathbf{P}, \mathbf{Q})+(2n-1)h(\mathbf{P}, \mathbf{Q})+\boldsymbol{\eta}(\mathbf{P})\boldsymbol{\eta}(\mathbf{Q}).
 \end{equation}
 Additionally, we substitute $\mathbf{P} = \tilde{e_i}$ and $\mathbf{Q} = \tilde{e_i}$ into the aforementioned expression, for which $\tilde{e_i}$ represents a local mutually orthogonal basis, and total with $\mathfrak{i}$ ranging from $1$ to $2n+1$ to obtain
 \begin{equation}
 \mathcal{R}^\ast = \mathcal{R} + 4 n^2
 \end{equation}
 wherein $\mathcal{R}^\ast$ corresponds to the modified scalar curvature($\ast-scalar$) of $\mathbb{M}$.

\section{Kenmotsu Metrics Admitting $\ast - \boldsymbol{\kappa}$-Ricci--Bourguignon Almost Solitons ($\ast - \boldsymbol{\kappa}\mathrm{RBS}$)}\label{3}
\hspace{0.5cm}
 Assume $\mathbb{M}$ to be a Kenmotsu-type manifold of $(2n+1)$-dim space. We now examine $\mathcal{V}=\zeta$ in identity (1.9) on $\mathbb{M}$ to obtain
 \begin{equation}
 \boldsymbol{\kappa}(\mathscr{L}_\zeta h)(\mathbf{P}, \mathbf{Q})+2\mathbb{T}^\ast(\mathbf{P}, \mathbf{Q})-2(\Omega+\theta \mathcal{R}^\ast)h(\mathbf{P}, \mathbf{Q})=0
 \end{equation}
  concerning any arbitrary vector fields $\mathbf{P}, \mathbf{Q} \in \mathfrak{X}(\mathbb{M})$.\\
Thus,the principal theorems are stated as follows
\begin{thm}\label{Th1}
   If the given tensor metric $h$ for a $(2n+1)$ odd dim-space Kenmotsu-type Manifold admits $\ast$- $\boldsymbol{\kappa}$-Ricci Bourguignon soliton($\ast-\boldsymbol{\kappa} RBS$) $(h,\zeta, \Omega,\theta)$, in which $\zeta$ equals canonical reeb vector field, it follows that geometric soliton becomes contracting, steady and  growing depending on whether $\theta(\mathcal{R}+4n^2) \gtreqqless 0$.
   \end{thm}
   \textbf{\textit{Proof.}}  By applying equations (2.14) and (2.15), equation (3.1) yields
   \begin{equation}
   \mathbb{T} (\mathbf{P}, \mathbf{Q})+[\Omega+(2n-\boldsymbol{\kappa})+1-\frac{\theta \mathcal{R}^\ast}{2}]h(\mathbf{P}, \mathbf{Q})+(1-\boldsymbol{\kappa})\boldsymbol{\eta}(\mathbf{P})\boldsymbol{\eta}(\mathbf{Q})=0
   \end{equation}
   Next, we substitute $\mathbf{Q} = \zeta$ into the before mentioned equation and using identities (2.1) and (2.10) to get
   \begin{equation}
   [\Omega-\frac{\theta \mathcal{R}^\ast}{2}]\boldsymbol{\eta}(\mathbf{P})=0
   \end{equation}
   Given that $\boldsymbol{\eta}(\mathbf{P}) \neq 0$, the preceding equation assumes the form of
   \begin{equation}
   \Omega = \frac{\theta \mathcal{R}^\ast}{2}
   \end{equation}
   Using (2.16), we obtain
   \begin{equation}
   \Omega = \frac{\theta (\mathcal{R}+4n^2)}{2}.
   \end{equation}
   Thus, we complete the proof..\\
   
   Based on the proven assertion, one deduces that
   
\begin{cor}
Let the tensor metric $h$ of a $(2n+1)$ odd-dim space Kenmotsu-type manifold, namely Riemann-flat, varifies the $\ast-\boldsymbol{\kappa} RBS$ $(h, \zeta, \Omega, \theta)$, as $\zeta$ denotes canonical reeb vector field, thus the geometric soliton is contracting, steady and growing with respect to $\theta \gtreqqless 0$.
\end{cor} 
When the $\mathbb{M}$ is Riemann-flat, which is $\mathcal{R}=0$, equation (3.5) simplifies to $\Omega = 2\theta n^2$, yielding the outcome.\\

Thereafter, we derive
\begin{thm}
If the tensor metric $h$ of a $(2n+1)$ odd-dim Kenmotsu-type manifold is endowed with the $\ast-\boldsymbol{\kappa} RBS$ $(h,\mathcal{V},\Omega,\theta)$, such that the function $\mathcal{V}$ represents gradient vector corresponding to continuously differentiable mapping $u$, therefore the Poisson's relation holds for $u$ signifies \\

$$\bigtriangleup (u)=- \frac{(\mathcal{R}+4n^2)}{\boldsymbol{\kappa}}[1+\frac{\theta(2n+1)}{2}]-\frac{\Omega}{\boldsymbol{\kappa}}(2n+1). $$
\end{thm}
\textbf{\textit{Proof.}} This work defines a $\ast$-$\boldsymbol{\kappa}$-Ricci-Bourguignon almost soliton($\ast-\boldsymbol{\kappa} RBS$) $(h, \mathcal{V}, \Omega, \theta)$ upon $\mathbb{M}$ in the manner follows:
\begin{equation}
 \boldsymbol{\kappa}(\mathscr{L}_\mathcal{V} h)(\mathbf{P}, \mathbf{Q})+2\mathbb{T}^\ast(\mathbf{P}, \mathbf{Q})-2(\Omega+\theta \mathcal{R}^\ast)h(\mathbf{P}, \mathbf{Q})=0
 \end{equation}
 for any vector fields $\mathbf{P}, \mathbf{Q}, \in \mathfrak{X} (\mathbb{M})$.\\
 Let $\mathbf{P} = \tilde{e_i}$ and $\mathbf{Q} = \tilde{e_i}$ within the preceding expression, subject to $\tilde{e_i}'s$ constitute a local mutually orthogonal basis. Following this, add up $\mathfrak{i} = 1, 2, ..., (2n + 1)$ and apply (2.16) to get
\begin{equation}
\mathcal{D}(\mathcal{V})+ \frac{\mathcal{R}+4n^2}{\boldsymbol{\kappa}}[1+\frac{\theta(2n+1)}{2}]-\frac{\Omega}{\boldsymbol{\kappa}}(2n+1)=0.
 \end{equation}
 Consider the vector function $\mathcal{V}$ is of gradient vector type, that is $\mathcal{V}  = \nabla(u)$, where $u$ denotes continuously differentiable function lying on $\mathbb{M}$, therefore the identity (3.7) transforms to
 \begin{equation}
 \bigtriangleup (u)=- \frac{(\mathcal{R}+4n^2)}{\boldsymbol{\kappa}}[1+\frac{\theta(2n+1)}{2}]-\frac{\Omega}{\boldsymbol{\kappa}}(2n+1)
 \end{equation}
 where $\bigtriangleup (u)$ denotes the  Poisson's relation defined by $u$. \\Thus, the proof is complete.\\ 
 Proceeding further, let us take $\theta=0$, (1.9) degenerates into $\ast$-$\boldsymbol{\kappa}$-Ricci flow soliton and (3.8) emerges as $\bigtriangleup(u)=-\frac{(\mathcal{R}+4n^2)}{\boldsymbol{\kappa}}-\frac{\Omega}{\boldsymbol{\kappa}}(2n+1)$.\\
 Given that $\theta=2$, (1.9) simplifies to $\ast$-$\boldsymbol{\kappa}$-Bourguignon geometric soliton with (3.8) is written as, $\bigtriangleup(u)=\frac{1}{\boldsymbol{\kappa}}[\mathcal{R}+4n^2-\Omega](2n+1)$.\\
Additionally, $\theta=1$, (1.9) transforms into $\ast$-$\boldsymbol{\kappa}$-Einstein-type soliton and (3.8) yields $\bigtriangleup(u)=-\frac{(\mathcal{R}+4n^2)}{\boldsymbol{\kappa}}[1+\frac{(2n+1)}{2}]-\frac{\Omega}{\boldsymbol{\kappa}}(2n+1).$ 

\begin{defn}
 Let the vector-valued function $\mathcal{V}$ is classified as a conformal isometry-generating vector field if the subsequent equation is satisfied:
\begin{equation}
    (\mathscr{L}_\mathcal{V} h)(\mathbf{P}, \mathbf{Q})=2\Omega{h}(\mathbf{P}, \mathbf{Q})
\end{equation}
 At which $\Omega$ represents a certain mapping in terms of this parameters (scaling scalar).\\
 In case $\theta$ signifies non-invariant, said conformal isometry-generating vector field $\mathcal{V}$ is termed as Proper .
 Whenever $\Omega$ is invariant, $\mathcal{V}$ referred to as a homothetic conformal killing vector field; at the time, that invariant $\Omega$ turns non-trivial, $\mathcal{V}$ is classified as a Proper homothetic conformal killing vector field. On condition that, $\Omega = 0$ in the foregoing-mentioned relation, therefore $\mathcal{V}$ corresponds to a infinitesimal isometry vector field.

\end{defn}
  
\begin{thm}
    If the tensor metric $h$ of a $(2n+1)$ odd-dim Kenmotsu-type manifold endows this $\ast$-$\boldsymbol{\kappa} RBS$ $(h, \mathcal{V}, \Omega, \theta)$, under which field $\mathcal{V}$ constitutes a conformal isometry-generating vector field, hence, this differentiable manifold is rendered $\boldsymbol{\eta}$-Einstein, provided $\alpha \neq 0$.
\end{thm}
\textbf{\textit{Proof.}}  As $(h, \mathcal{V}, \Omega, \theta)$ constitutes one $\ast$-$\boldsymbol{\kappa}$-Ricci Bourguignon almost soliton on an $(2n+1)$ odd-dim Kenmotsu-type manifold $\mathbb{M}$, with $\mathcal{V}$ representing a conformal
isometry-generating vector field.
Consequently, based on (1.9), (2.15), and (3.9), we derive
\begin{equation}
    \mathbb{T}(\mathbf{P}, \mathbf{Q})=-[(2n-1)-\Omega+\boldsymbol{\kappa}\Lambda+\frac{\theta \mathcal{R}^\ast}{2}]h(\mathbf{P}, \mathbf{Q})-\boldsymbol{\eta}(\mathbf{P})\boldsymbol{\eta}(\mathbf{Q}).
\end{equation}
Thus it is understood, an manifold of Kenmotsu-type $(M^{2n+1}, h)$ is categorized as $\boldsymbol{\eta}$-Einstein subject to the respective tensor of Ricci curvature $\mathbb{T}$ having order (0, 2) conforms to such above expression.
\begin{equation}
    \mathbb{T}=\alpha h+\beta \boldsymbol{\eta} \otimes \boldsymbol{\eta}.
\end{equation}
In which $\alpha$ and $\beta$ are taken as continuosly differentiable  function lying on $(\mathbb{M}^{2n+1}, h)$.\\
Since (3.11), (3.10) determines with the property that the specified geometric manifold proves to be $\boldsymbol{\eta}$-Einstein, assuming that $\alpha \neq 0$, Thus, the statement holds.
\begin{thm}
    Suppose the tensor metric $h$ for a $(2n+1)$ odd-dim Kenmotsu-type manifold induces on the $\ast$-$\boldsymbol{\kappa} RBS$ $(h, \mathcal{V}, \Omega, \theta)$, such that $\mathcal{V}$ denotes the conformal isometry-generating vector field . Accordingly, $\mathcal{V}$ represents\\
\textbf{(i)} Proper tangent vector field only if $\frac{1}{\boldsymbol{\kappa}}[\frac{\theta(\mathcal{R}+4n^2)}{2}+\Omega$ denotes non-invariant.\\
 \textbf{(ii)} Homothetic conformal killing vector field if $\frac{1}{\boldsymbol{\kappa}}[\frac{\theta(\mathcal{R}+4n^2)}{2}+\Omega$ signifies invariant.\\
\textbf{(iii)} Proper Homothetic conformal killing tangent vector field only if $\frac{1}{\boldsymbol{\kappa}}[\frac{\theta(\mathcal{R}+4n^2)}{2}+\Omega$ represents non-trivial invariant.\\
\textbf{(iv)} Infinitesimal isometry vector field only if $\Omega=\frac{\theta(\mathcal{R}+4n^2)}{2}$.

\end{thm}
\textbf{\textit{Proof.}}  We substitute $\mathbf{Q} = \zeta$ into (3.10) and apply (2.1) and (2.10) to obtain
\begin{equation}
    [2n-(2n-1)+\Omega-\boldsymbol{\kappa} \Lambda-\frac{\theta \mathcal{R}^\ast}{2}-1]\boldsymbol{\eta}{\mathbf{P}}=0
\end{equation}
Given that $\boldsymbol{\eta}(\mathbf{P}) \neq 0$, we have
\begin{equation}
    \Lambda=\frac{1}{\boldsymbol{\kappa}}[\frac{\theta \mathcal{R}^\ast}{2}-\Omega]
\end{equation}
Following that, use equation (2.16), (3.13) can be stated as follows
\begin{equation}
    \Lambda=\frac{1}{\boldsymbol{\kappa}}[\frac{\theta (\mathcal{R}+4n^2)}{2}-\Omega].
\end{equation}
Hence, we get the conclusion of the result.
\begin{thm}
    Given that $(\mathbb{M}, \varphi, h, \zeta, \Omega, \theta, \alpha, \beta)$ is $\ast$-$\boldsymbol{\kappa}$-Ricci-Bourguignon almost soliton$(\ast-\boldsymbol{\kappa}RBS)$ concerning a certain $\boldsymbol{\eta}$-Einstein Kenmotsu-type manifold,it follows that $\Omega=-\alpha-\frac{\theta}{2}(\mathcal{R}+4n^2)-2n-\beta$
\end{thm}
\textbf{Proof.}  We combine identities (3.1), (2.15), and (2.16) to derive

\begin{equation}\label{3.15}
\begin{split}
\boldsymbol{\kappa}\big[h(\nabla_\mathbf{P} \mathcal{V}, \mathbf{Q}) + h(\mathbf{P}, \nabla_\mathcal{V} \mathbf{Q})\big] + 2\mathbb{T}(\mathbf{P}, \mathbf{Q}) + \big[2(2n - 1) \\
\hspace*{\fill} - 2\big[\Omega + \theta(\mathcal{R} + 4n^2)]] h(\mathbf{P}, \mathbf{Q}) + 2\boldsymbol{\eta}(\mathbf{P})\boldsymbol{\eta}(\mathbf{Q}) = 0
\end{split}
\end{equation}
We add relations (3.11) along with (3.10) within the above formula to find

\begin{equation}\label{3.16}
\begin{split}
\boldsymbol{\kappa}\big[h(\nabla_\mathbf{P} \mathcal{V}, \mathbf{Q}) + h(\mathbf{P}, \nabla_\mathcal{V} \mathbf{Q})\big] + [2 \alpha - 2\Omega + \theta(\mathcal{R} + 4n^2) \\
\hspace*{\fill} + 2(2n - 1)]h(\mathbf{P}, \mathbf{Q}) + [2(\beta + 1)]\boldsymbol{\eta}(\mathbf{P})\boldsymbol{\eta}(\mathbf{Q}) = 0
\end{split}
\end{equation}
We now put $\mathbf{P} = \mathbf{Q} = \zeta$ into equation (3.16) so that\\
\begin{equation}
2\boldsymbol{\kappa} h(\nabla_\zeta \mathcal{V}, \zeta)=[2 \alpha-2\Omega+\theta(\mathcal{R}+4n^2)+2(2n-1)]+2(\beta+1) 
\end{equation}
Let place $\mathcal{V} = \zeta$ within the relation (3.17) so as to generate\\
\begin{equation}
\boldsymbol{\kappa} h(\nabla_\zeta \zeta, \zeta)=[\alpha-\Omega+\frac{\theta}{2}(\mathcal{R}+4n^2)+2n+\beta]
\end{equation}
Since it is commonly known that\\
$$h(\nabla_\zeta \zeta, \zeta)=0$$\\
Regarding each tangent fields associated with the manifold $\mathbb{M}$, it is deduced that\\
$$\Omega=-\alpha-\frac{\theta}{2}(\mathcal{R}+4n^2)-2n-\beta.$$\\
Thus, the Proof.

\section{On the Role arising from Torse-Generating tangent Fields in $\ast$-$\boldsymbol{\kappa}$-Ricci–Bourguignon Almost Soliton on Kenmotsu-type Manifolds}\label{4}
A topic of this article will focus on the role of torsion-generating potential tangent fields in $\ast$-$\boldsymbol{\kappa}RBS$ on Kenmotsu-type manifolds. Herein, we provide the succeeding theorem.
\begin{thm}
    If the tensor metric $h$ related to  a $(2n+1)$ odd-dim Kenmotsu-
type manifold possesses the $\ast$-$\boldsymbol{\kappa}RBS$ $(h, \rho, \Omega, \theta)$, where $\rho$ is a torse-forming vector field , then $\Omega= \psi \boldsymbol{\kappa}- (2n-1)-\frac{\theta}{2}(\mathcal{R}+4n^2)-\frac{\mathcal{R}+1+\boldsymbol{\kappa}\theta(\rho)}{2n+1}$ and the soliton is contracting, steady, growing according to $\psi \boldsymbol{\kappa}- (2n-1)-\frac{\theta}{2}(\mathcal{R}+4n^2)-\frac{\mathcal{R}+1+\boldsymbol{\kappa}\theta(\rho)}{2n+1} \gtreqqless 0.$
\end{thm}
\textbf{\textit{Proof.}}  Let $(h, \rho, \Omega, \theta)$ represent an $\ast$-$\boldsymbol{\kappa}$-Ricci Bourguignon almost soliton($\ast$-$\boldsymbol{\kappa}RBS$) on a $(2n + 1)$odd-dim Kenmotsu-
type $\mathbb{M}$, with $\rho$ providing as a torse-generating vector field.\\
 Consequently, from equations (1.9), (2.15), and (2.16), we obtain
\begin{equation}
\begin{split}
\boldsymbol{\kappa}(\mathscr{L}_\rho h)(\mathbf{P}, \mathbf{Q})+2[\mathbb{T}(\mathbf{P}, \mathbf{Q})+(2n-1)h(\mathbf{P}, \mathbf{Q})+\boldsymbol{\eta}(\mathbf{P})\boldsymbol{\eta}(\mathbf{Q})] \\ -[2\Omega+\theta(\mathcal{R}+4n^2)]h(\mathbf{P}, \mathbf{Q})=0
\end{split}
\end{equation}
In case where $\mathscr{L}_\rho h$ indicates the Lie derivative of the metric $h$ overlying the tangent field $\rho$.\\
We now employ equation (1.10) to determine
\begin{multline} \label{4.2}
(\mathscr{L}_\rho h)(\mathbf{P}, \mathbf{Q}) = h(\nabla_\mathbf{P} \rho, \mathbf{Q}) + h(\mathbf{P}, \nabla_\mathbf{Q} \rho) \\
\hfill = 2\psi\, h(\mathbf{P}, \mathbf{Q}) + \theta(\mathbf{P})\, h(\rho, \mathbf{Q}) + \theta(\mathbf{Q})\, h(\rho, \mathbf{P})
\end{multline}
for all $\mathbf{P}, \mathbf{Q} \in \mathbb{M}$.\\
We further combine identities (4.3) and (4.1) to form
\begin{equation}\label{4.3}
\begin{split}
\left[\psi \boldsymbol{\kappa} - (2n - 1)+ \Omega-\frac{\theta(\mathcal{R} + 4n^2)}{2}\right] h(\mathbf{P}, \mathbf{Q})
- \mathbb{T}(\mathbf{P}, \mathbf{Q}) - \boldsymbol{\eta}(\mathbf{P})\boldsymbol{\eta}(\mathbf{Q}) \\
\hspace*{\fill} =\frac{\boldsymbol{\kappa}}{2} \left[\theta(\mathbf{P}) h(\rho, \mathbf{Q}) + \theta(\mathbf{Q}) h(\rho, \mathbf{P})\right]
\end{split}
\end{equation}
We now do the contraction of (4.3) over $\mathbf{P}$ and $\mathbf{Q}$ so that
\begin{equation}
[\psi \boldsymbol{\kappa} - (2n - 1)+ \Omega-\frac{\theta(\mathcal{R} + 4n^2)}{2}](2n+1)-\mathcal{R}-1=\boldsymbol{\kappa}\theta(\rho)
\end{equation}
Which ends up in
\begin{equation}
\Omega=\psi \boldsymbol{\kappa}- (2n-1)-\frac{\theta(\mathcal{R}+4n^2)}{2}-\frac{\mathcal{R}+1+\theta(\rho)}{2n+1}.
\end{equation}
Based on the characteristics with respect to the soliton, this work establish the last principle concerning the theorem.\\
Thus, the Proof is completed.\\
   As seen in relation (4.5), when this linear-functional(1-form) form $\theta$ is exactly zero, thus $\Omega=\psi \boldsymbol{\kappa}-(2n-1)-\frac{(\mathcal{R}+1)}{2n+1}$.\\
   Assuming that, (1-form) $\theta$ is exactly zero also the given mapping $\psi = 1$ within (4.5), and hence $\Omega=\boldsymbol{\kappa}-(2n-1)-\frac{\mathcal{R}+1}{2n+1}$.\\
   If the function $\psi = 0$ in (4.5), then $\Omega=(2n-1)-\frac{\theta(\mathcal{R}+4n^2)}{2}-\frac{\mathcal{R}+1+\boldsymbol{\kappa}\theta(\rho)}{2n+1}$.\\
   In the case where $\psi = \theta = 0$ in (4.5), $\Omega=(2n-1)-\frac{\mathcal{R}+1}{2n+1}$.\\
   To finalize (4.5), $\theta(\rho)=0$, then $\Omega=\psi \boldsymbol{\kappa}-(2n-1)-\frac{\theta(\mathcal{R}+4n^2)}{2}-\frac{\mathcal{R}+1}{2n+1}$.\\
   Furthermore, we include
   \begin{cor}
       Consider the tensor metric $h$ associated with  a $(2n+1)$ odd-dim Kenmotsu-type manifold  fulfills the $\ast$-$\boldsymbol{\kappa}RBS$ $(h, \rho, \Omega, \theta)$, where $\rho$ acts as a torse-generating tangent field, therefore if $\rho$ can be identified as \\
       
       \textbf{(i)}conformally circular, thus $\Omega=\psi \boldsymbol{\kappa}-(2n-1)-\frac{\mathcal{R}+1}{2n+1}$ and the soliton is contracting, steady, growing according as $\Omega=\psi \boldsymbol{\kappa}-(2n-1)-\frac{\mathcal{R}+1}{2n+1}\gtreqqless 0,$\\
       
       \textbf{(ii)} concurrent,then $\Omega=\boldsymbol{\kappa}-(2n-1)-\frac{\mathcal{R}+1}{2n+1}$ and the soliton is contracting, steady, growing according as $\Omega=\boldsymbol{\kappa}-(2n-1)-\frac{\mathcal{R}+1}{2n+1}\gtreqqless 0,$\\
       
       \textbf{(iii)}recurrent, then $\Omega=(2n-1)-\frac{\theta(\mathcal{R}+4n^2)}{2}-\frac{\mathcal{R}+1+\boldsymbol{\kappa}\theta(\rho)}{2n+1}$ and the soliton is contracting, steady, growing according as $\Omega=(2n-1)-\frac{\theta(\mathcal{R}+4n^2)}{2}-\frac{\mathcal{R}+1+\boldsymbol{\kappa}\theta(\rho)}{2n+1}\gtreqqless 0$,\\
       
       \textbf{(iv)}covariantly constant, therefore $\Omega=(2n-1)-\frac{\mathcal{R}+1}{2n+1}$ also the specified geometric soliton represents contracting, steady, growing according as $\Omega=(2n-1)-\frac{\mathcal{R}+1}{2n+1}\gtreqqless 0$, \\

       \textbf{(v)}torqued, then $\Omega=\psi \boldsymbol{\kappa}-(2n-1)-\frac{\theta(\mathcal{R}+4n^2)}{2}-\frac{\mathcal{R}+1}{2n+1}$ and the soliton is contracting, steady, growing according as $\Omega=\psi \boldsymbol{\kappa}-(2n-1)-\frac{\theta(\mathcal{R}+4n^2)}{2}-\frac{\mathcal{R}+1}{2n+1}\gtreqqless 0.$
   \end{cor}

\section{"Explicit Illustration belonging to a 5D Kenmotsu Structure Satisfying the $\ast$-$\boldsymbol{\kappa}$-Ricci-Bourguignon Almost Soliton Equation"}\label{5}

\textbf{Illustration 6.1.}
    Introduce the set $\mathbb{M} = (\mathbf{x_1, x_2, x_3, x_4, x_5}) \in \mathcal{R}^5$ taken as our manifold, such that $(\mathbf{x_1, x_2, x_3, x_4, x_5})$ denote the ordinary coordinates embedded in $\mathcal{R}^5$.  The tangent fields are shown in the subsequent:
    \[
\begin{aligned}
\tilde{e_1} &= \tilde{e}^{-\mathbf{x_5}} \frac{\partial}{\partial \mathbf{x_1}}, \quad
\tilde{e_2} = \tilde{e}^{-\mathbf{x_5}} \frac{\partial}{\partial \mathbf{x_2}}, \quad
\tilde{e_3} = \tilde{e}^{-\mathbf{x_5}} \frac{\partial}{\partial \mathbf{x_3}}, \quad
\tilde{e_4} = \tilde{e}^{-\mathbf{x_5}} \frac{\partial}{\partial \mathbf{x_4}}, \quad
\tilde{e_5} = \frac{\partial}{\partial \mathbf{x_5}}
\end{aligned}
\]
constitute a mutually independent system at every element of $\mathbb{M}$.  One determine the tensor metric h given by
\[
h(\tilde{e_i}, \tilde{e_j}) = 
\begin{cases}
1, & \text{in the case that } i = j \text{ and } i, j \in \{1, 2, 3, 4, 5\} \\
0, & \text{in any other case}
\end{cases}
\]
Take $\boldsymbol{\eta}$ constitute a linear-functional(1-form) form expressed as $\boldsymbol{\eta}(\mathbf{P}) = h(\mathbf{P}, \tilde{e_5})$ in order to every $\mathbf{P} \in \mathfrak{X}(\mathbb{M})$.  We describe the (1,1)-tensorial field $\varphi$ in the form of:
\[
\begin{aligned}
\varphi\tilde{(e_1)} = \tilde{e_3}, \quad  \varphi\tilde{(e_2)} = \tilde{e_4}, \quad
\varphi\tilde{(e_3)} = \tilde{-e_1}, \quad  \varphi\tilde{(e_4)} = \tilde{-e_2}, \quad \varphi\tilde{(e_5)} = 0.
\end{aligned}
\]
Thus, this satisfies the conditions $\boldsymbol{\eta}(\zeta)=1$, $\varphi^2(\mathbf{P})=\mathbf{-P}+\boldsymbol{\eta}(\mathbf{P})\zeta$   also  $h(\varphi \mathbf{P}, \varphi \mathbf{Q})=h(\mathbf{P}, \mathbf{Q})-\boldsymbol{\eta}(\mathbf{P})\boldsymbol{\eta}(\mathbf{Q})$ , for which $\zeta = \tilde{e_5}$ and $\mathbf{P}, \mathbf{Q}$ are any vector fields over $\mathbb{M}$. $(\mathbb{M}, \varphi, \zeta, \boldsymbol{\eta}, h)$ constitutes a certain almost contact geometric structure in $\mathbb{M}$.\\
\noindent\hspace*{0.5cm}It may at this stage be inferred which
\[
\begin{aligned}
\relax [\tilde{e_1}, \tilde{e_2}] &= 0, & [\tilde{e_1}, \tilde{e_3}] &= 0, & [\tilde{e_1}, \tilde{e_4}] &= 0, & [\tilde{e_1}, \tilde{e_5}] &= \tilde{e_1}, \\
[\tilde{e_2}, \tilde{e_1}] &= 0, & [\tilde{e_2}, \tilde{e_3}] &= 0, & [\tilde{e_2}, \tilde{e_4}] &= 0, & [\tilde{e_2}, \tilde{e_5}] &= \tilde{e_2}, \\
[\tilde{e_3}, \tilde{e_1}] &= 0, & [\tilde{e_3}, \tilde{e_2}] &= 0, & [\tilde{e_3}, \tilde{e_4}] &= 0, & [\tilde{e_3}, \tilde{e_5}] &= \tilde{e_3}, \\
[\tilde{e_4}, \tilde{e_1}] &= 0, & [\tilde{e_4}, \tilde{e_2}] &= 0, & [\tilde{e_4}, \tilde{e_3}] &= 0, & [\tilde{e_4}, \tilde{e_5}] &= \tilde{e_4}, \\
[\tilde{e_5}, \tilde{e_1}] &= \tilde{-e_1}, & [\tilde{e_5}, \tilde{e_2}] &= \tilde{-e_2}, & [\tilde{e_5}, \tilde{e_3}] &= \tilde{-e_3}, & [\tilde{e_5}, \tilde{e_4}] &= \tilde{-e_4}.
\end{aligned}
\]
Denote $\nabla$ show the particular Levi-Civita connection associated with $h$.  Subsequently, from Koszul's formula for arbitrary $\mathbf{P}, \mathbf{Q}, \mathbf{R} \in \mathfrak{X}(\mathbb{M})$ as provided:
\[
\begin{aligned}
2h(\nabla_\mathbf{P} \mathbf{Q}, \mathbf{R}) =\ & \mathbf{P} h(\mathbf{Q},\mathbf{R}) + \mathbf{Q} h(\mathbf{R},\mathbf{P}) - \mathbf{R}h(\mathbf{P},\mathbf{Q}) - h(\mathbf{P}, [\mathbf{Q},\mathbf{R}]) \\
& - h(\mathbf{Q}, [\mathbf{P},\mathbf{R}]) + h(\mathbf{R}, [\mathbf{P},\mathbf{Q}]).
\end{aligned}
\]
It may be already establish:
\[
\begin{aligned}
\nabla_{\tilde{e_1}} \tilde{e_1} &= \tilde{-e_5}, & \nabla_{\tilde{e_1}} \tilde{e_2} &= 0, & \nabla_{\tilde{e_1}} \tilde{e_3} &= 0, & \nabla_{\tilde{e_1}} \tilde{e_4} &= 0, & \nabla_{\tilde{e_1}} \tilde{e_5} &= \tilde{e_1}, \\
\nabla_{\tilde{e_2}} \tilde{e_1} &= 0,   & \nabla_{\tilde{e_2}} \tilde{e_2} &= \tilde{-e_5}, & \nabla_{\tilde{e_2}} \tilde{e_3} &= 0, & \nabla_{\tilde{e_2}} \tilde{e_4} &= 0, & \nabla_{\tilde{e_2}} \tilde{e_5} &= \tilde{e_2}, \\
\nabla_{\tilde{e_3}} \tilde{e_1} &= 0,   & \nabla_{\tilde{e_3}} \tilde{e_2} &= 0,   & \nabla_{\tilde{e_3}} \tilde{e_3} &= \tilde{-e_5}, & \nabla_{\tilde{e_3}} \tilde{e_4} &= 0, & \nabla_{\tilde{e_3}} \tilde{e_5} &= \tilde{e_3}, \\
\nabla_{\tilde{e_4}} \tilde{e_1} &= 0,   & \nabla_{\tilde{e_4}} \tilde{e_2} &= 0,   & \nabla_{\tilde{e_4}} \tilde{e_3} &= 0, & \nabla_{\tilde{e_4}} \tilde{e_4} &= \tilde{-e_5}, & \nabla_{\tilde{e_4}} \tilde{e_5} &= \tilde{e_4}, \\
\nabla_{\tilde{e_5}} \tilde{e_1} &= 0,   & \nabla_{\tilde{e_5}} \tilde{e_2} &= 0,   & \nabla_{\tilde{e_5}} \tilde{e_3} &= 0, & \nabla_{\tilde{e_5}} \tilde{e_4} &= 0, & \nabla_{\tilde{e_5}} \tilde{e_5} &= 0.
\end{aligned}
\]
Consequently, $(\nabla_\mathbf{P}\varphi)\mathbf{Q} = h(\varphi \mathbf{P}, \mathbf{Q})\zeta - \boldsymbol{\eta}(\mathbf{Q})\varphi \mathbf{P}$ holds for every $\mathbf{P}, \mathbf{Q} \in \mathfrak{X}(\mathbb{M})$.  $(\mathbb{M}, \varphi, \zeta, \boldsymbol{\eta}, h)$ acts as a Kenmotsu-type manifold.\\
These non-trivial coordinates of the Riemannian curvature tensor takes the form:
\[
\begin{aligned}
\mathcal{R}(\tilde{e_1}, \tilde{e_2})\tilde{e_2} &= \tilde{-e_1}, & \mathcal{R}(\tilde{e_1}, \tilde{e_3})\tilde{e_3} &= \tilde{-e_1}, & \mathcal{R}(\tilde{e_1}, \tilde{e_4})\tilde{e_4} &= \tilde{-e_1}, \\
\mathcal{R}(\tilde{e_1}, \tilde{e_5})\tilde{e_5} &= \tilde{-e_1}, & \mathcal{R}(\tilde{e_1}, \tilde{e_2})\tilde{e_1} &= \tilde{e_2}, & \mathcal{R}(\tilde{e_1}, \tilde{e_3})\tilde{e_1} &= \tilde{e_3}, \\
\mathcal{R}(\tilde{e_1}, \tilde{e_4})\tilde{e_1} &= \tilde{e_4}, & \mathcal{R}(\tilde{e_1}, \tilde{e_5})\tilde{e_1} &= \tilde{e_5}, & \mathcal{R}(\tilde{e_2}, \tilde{e_3})\tilde{e_2} &= \tilde{e_3}, \\
\mathcal{R}(\tilde{e_2}, \tilde{e_4})\tilde{e_2} &= \tilde{e_4}, & \mathcal{R}(\tilde{e_2}, \tilde{e_5})\tilde{e_2} &= \tilde{e_5}, & \mathcal{R}(\tilde{e_2}, \tilde{e_3})\tilde{e_3} &= \tilde{-e_2}, \\
\mathcal{R}(\tilde{e_2}, \tilde{e_4})\tilde{e_4} &= \tilde{-e_2}, & \mathcal{R}(\tilde{e_2}, \tilde{e_5})\tilde{e_5} &= \tilde{-e_2}, & \mathcal{R}(\tilde{e_3}, \tilde{e_4})\tilde{e_3} &= \tilde{e_4}, \\
\mathcal{R}(\tilde{e_3}, \tilde{e_5})\tilde{e_3} &= \tilde{e_5}, & \mathcal{R}(\tilde{e_3}, \tilde{e_4})\tilde{e_4} &= \tilde{-e_3}, & \mathcal{R}(\tilde{e_4}, \tilde{e_5})\tilde{e_4} &= \tilde{e_5}, \\
\mathcal{R}(\tilde{e_5}, \tilde{e_3})\tilde{e_5} &= \tilde{e_3}, & \mathcal{R}(\tilde{e_5}, \tilde{e_4})\tilde{e_5} &= \tilde{e_4}. & &
\end{aligned}
\]
On the basis of such prior outcomes, it can get $\mathbb{T}\tilde{(e_i}, \tilde{e_i}) = -4$ with respect to $i = 1, 2, 3, 4, 5$.
\begin{equation}
\mathbb{T}(\mathbf{P}, \mathbf{Q}) = -4h(\mathbf{P}, \mathbf{Q}) \quad \forall\, \mathbf{P}, \mathbf{Q} \in \mathfrak{X}(\mathbb{M}).
\end{equation}
By contracting, we get  $\mathcal{R} = \sum_{i=1}^{5}$ $\mathbb{T}(\tilde{e_i}, \tilde{e_i}) = -20 = -2n(2n+1)$, such that the dimension belonging to the present manifold is $2n + 1 = 5$. Moreover, one know
\[
\mathbb{T}(\tilde{e_i}, \tilde{e_i}) = 
\begin{cases}
-1, & \text{in case } $i = 1, 2, 3, 4$ \\
0, & \text{in case } $i = 5$
\end{cases}
\]
Thus,
\begin{equation}
\mathbb{T}^\ast(\mathbf{P}, \mathbf{Q}) = -h(\mathbf{P}, \mathbf{Q}) + \boldsymbol{\eta}(\mathbf{P})\boldsymbol{\eta}(\mathbf{Q}) \quad \forall\, \mathbf{P}, \mathbf{Q} \in \mathfrak{X}(\mathbb{M}).
\end{equation}
Since
\begin{equation}
\mathcal{R}^\ast = \operatorname{Tr}(\mathbb{T}^\ast) = -4,
\end{equation}
We will now examine a tangent field $\mathcal{V}$ known by the term
\begin{equation}
\mathcal{V} = \mathbf{x_1} \frac{\partial}{\mathbf{x_1}} + \mathbf{x_2} \frac{\partial}{\mathbf{x_2}} + \mathbf{x_3} \frac{\partial}{\mathbf{x_3}} + \mathbf{x_4} \frac{\partial}{\mathbf{x_4}} + \frac{\partial}{\mathbf{x_5}}
\end{equation}
Therefore, according to the results above, we can confirm that
\begin{equation}
(\mathscr{L}_\mathcal{V} h)(\mathbf{P}, \mathbf{Q}) = 4 \left\{ h(\mathbf{P}, \mathbf{Q}) - \boldsymbol{\eta}(\mathbf{P}) \boldsymbol{\eta}(\mathbf{Q}) \right\}.
\end{equation}
It's related for every $\mathbf{P}, \mathbf{Q} \in \mathfrak{X}(\mathbb{M})$.  Therefore, with respect to (5.5), we assume
\begin{equation}
\sum_{i=1}^{5} \bigl(\mathscr{L}_\mathcal{V} h\bigr)(\tilde{e_i}, \tilde{e_i}) \;=\; 16,
\end{equation}
We now set $\mathbf{P} = \mathbf{Q} = \tilde{e_i}$ in equation (1.9), summation on i = 1, 2, 3, 4, 5, and using equations (5.3) and (5.6) to obtain
\begin{equation}
    \Omega=\frac{4+10\theta-8\boldsymbol{\kappa}}{5}.
\end{equation}
Given that this $\Omega$, as stated above, fulfills (3.7), $h$ constitutes a $\ast$-$\boldsymbol{\kappa}$-Ricci-Bourguignon almost soliton($\ast$-$\boldsymbol{\kappa}RBS$) to a 5D Kenmotsu structure manifold $\mathbb{M}$.\\
Furthermore, it may be formulated as\\
\textbf{Remark 6.2.}\textbf{Case-I:} Whenever $\theta = 0$, equation (5.7) yields $\Omega = \frac{4 - 8\boldsymbol{\kappa}}{5}$, showing that $(h, \mathcal{V}, \Omega)$ makes a $\ast$-$\boldsymbol{\kappa}$-Ricci soliton, which is contracting for $\boldsymbol{\kappa} > \frac{1}{2}$, growing for $\boldsymbol{\kappa} < \frac{1}{2}$, and steady when $\boldsymbol{\kappa} = \frac{1}{2}$.\\
\textbf{Case-II:} When $\theta = 2$, equation (5.7) yields $\Omega = \frac{24 - 8\boldsymbol{\kappa}}{5}$, showing that $(h, \mathcal{V}, \Omega)$ makes a $\ast$-$\boldsymbol{\kappa}$-Bourguignon soliton, i.e., contracting provided that $\boldsymbol{\kappa} > \frac{-5}{2}$, growing whenever $\boldsymbol{\kappa} < \frac{-5}{2}$, as well as steady in case $\boldsymbol{\kappa} = \frac{-5}{2}$.\\
\textbf{Case-III:} If $\theta = 1$, equation (5.7) yields $\Omega=\frac{14-8\boldsymbol{\kappa}}{5}$, showing that the triplet $(h, \mathcal{V}, \Omega)$ makes a $\ast$-$\boldsymbol{\kappa}$-Einstein soliton, which is contracting for $\boldsymbol{\kappa}>-\frac{3}{4}$, growing for $\boldsymbol{\kappa} < -\frac{3}{4}$, and steady for $\boldsymbol{\kappa} = -\frac{3}{4}$.
 
\section*{Author Contributions}
\noindent
All authors contributed equally.

\section*{Data Availability Statement}
\noindent
Data sharing does not apply to this article, as no new data were generated or analyzed in this study.


\section*{Use of Generative-AI tools declaration}
\noindent
The authors declare they have not used artificial intelligence (AI) tools in the creation of this article.

\section*{Conflicts of Interest}
\noindent
The authors declare no conflicts of interest.

\section*{Funding} No external funding received for this research.

\noindent\textit{Corresponding author: Soumendu Roy.}

\hspace{0.5cm}



\begin{thebibliography}{99}

\bibitem{BYC}B.Y. Chen.: {\it A simple characterization of generalized Robertson-Walker space-times}, Gen. Relativity Gravitation {\bf 46} (12), Article No: 1833, 2014.

 \bibitem{BYC1}B.Y. Chen.: {\it Classification of torqued vector fields and its applications to Ricci solitons}, Kragujevac J. Math. {\bf 41} (2), 239-250, 2017.

\bibitem{BPS}C.S. Bagewadi and V. S. Prasad.: {\it Note on Kenmotsu manifolds}, Bull. Cal. Math. Soc. {\bf 91}, 379-384, 1999.

\bibitem{CDM}Catino, G., Cremaschi, L., Djadli, Z., Mantezza, C., Mazzieri, L.: {\it The Ricci- Bourguignon flow}. Pac. J. Math. {\bf 287}, 337370 (2017).

\bibitem{DPS}Dwibedi, S., Patra, D.S.: {\it Some results on almost $\ast$-Ricci-Bourguignon solitons}. J. Geom. Phys. {\bf 178}, 104519 (2022) 

\bibitem{DS}Dwivedi, S.: {\it Some results on Ricci-Bourguignon solitons and almost solitons}, (2020) {\bf arXiv:1809.11103v2} [math.DG]

\bibitem{PKG}G. Kaimakamis and K. Panagiotidou.: {\it *-Ricci Solitons of real hypersurface in non-flat complex space forms}, J. Geom. Phys. {\bf 76}, 408-413, 2014.

\bibitem{SJ}J.A. Schouten.: {\it Ricci Calculus}, Springer-Verlag, Berlin, 1954.

\bibitem{KK}K. Kenmotsu.: {\it A class of almost contact Riemannian manifolds}, The Tˆohoku Mathematical Journal, {\bf 24}, 93–103, 1972.

\bibitem{CYK}K. Yano and B.Y. Chen.: {\it On the concurrent vector fields of immersed manifolds}, Kodai Math. Sem. Rep. {\bf 23}, 343-350, 1971.

\bibitem{YK}K. Yano.: {\it Concircular geometry I. Concircular transformations}, Proc. Imp. Acad. Tokyo., {\bf 16}, 195-200, 1940.

\bibitem{YK1}K. Yano.: {\it On the torse-forming directions in Riemannian spaces}, Proc. Imp. Acad. Tokyo., {\bf 20}, 340345, 1944.

\bibitem{TP}P. Topping.: {\it Lecture on the Ricci Flow}, Cambridge University Press, 2006.

\bibitem{XGW}Q. Wang, J.N. Gomes and C. Xia.: {\it On the h-almost Ricci soliton}, J. Geom. Phys. {\bf 114}, 216-222, 2017.

\bibitem{HRS}R.S. Hamilton.: {\it The Ricci flow on surfaces}, Contemp. Math. {\bf 71}, 237-261, 1988.

\bibitem{HRS1}R.S. Hamilton.: {\it Three Manifold with positive Ricci curvature}, J. Differential Geom. {\bf 17} (2), 255-306, 1982.

\bibitem{BRS}S. Roy and A. Bhattacharyya.: {\it Conformal Ricci solitons on 3-dimensional trans-Sasakian manifold}, Jordan J. Math. Stat. {\bf 13} (1), 89-109, 2020.

\bibitem{BDR}S. Roy, S. Dey and A. Bhattacharyya.: {\it A Kenmotsu metric as a conformal-Einstein soliton}, Carpathian Math. Publ. {\bf 13} (1), 110-118, 2021.

\bibitem{TS}S. Tachibana.: {\it On almost-analytic vectors in almost KÅNahlerian manifolds}, Tohoku Math. J. {\bf 11} (2), 247-265, 1959.

\bibitem{TH}T. Hamada.: {\it Real hypersurfaces of complex space forms in terms of Ricci -tensor}, Tokyo J. Math. {\bf 25}, 473-483, 2002.

\bibitem{KNV}V. Venkatesha, D.M. Naik and H.A. Kumara.: {\it *-Ricci solitons and gradient almost-Ricci solitons on Kenmotsu manifolds}, {\bf arXiv:1901.05222} [math.DG].

\end{thebibliography}
\end{document}